\documentclass[12pt]{amsart}

\usepackage{amssymb}
\usepackage{epsf}
\usepackage{amsmath,epsf}
\usepackage[latin1]{inputenc}
\usepackage{amsfonts}

\numberwithin{equation}{section}
\newtheorem{thm}{Theorem}[section]
\newtheorem{prop}[thm]{Proposition}

\newtheorem{cor}[thm]{Corollary}

\newtheorem{rem}[thm]{Remark}

\def\Q{{\mathbb Q}}

\def\S{{\mathcal S}}

\def\Young#1{\vbox{\smallskip\offinterlineskip
    \halign{&\vbox{##}\kern-\Thickness\cr #1}}}

\newdimen\Squaresize \Squaresize=20pt
\newdimen\Thickness \Thickness=.1pt
\newdimen\Correction \Correction=7pt

\def\Vide#1{\hbox{
       \vbox to \Squaresize{\vss
          \hbox to \Squaresize{\hss#1 \hss}\vss}
    \hskip-\Correction}
   \kern-\Thickness}

\def\Carre#1{\hbox{\vrule width \Thickness
   \vbox to \Squaresize{\hrule height \Thickness\vss
      \hbox to \Squaresize{\hss$\scriptstyle#1$\hss}
   \vss\hrule height\Thickness}
   \unskip\vrule width \Thickness}
   \kern-\Thickness}

\title[Schur Operator]{Schur Partial Derivative Operators}

\author[J.-C.~Aval and N.~Bergeron]
{J.-C.~Aval and N.~Bergeron}

\address[Jean-Christophe Aval]
{Laboratoire A2X\\ Universit\'e Bordeaux 1\\ 351 cours de la
Lib\'eration\\ 33405 Talence cedex\\ FRANCE}

\address[Nantel Bergeron]
{Department of Mathematics and Statistics\\ York University\\
   To\-ron\-to, Ontario M3J 1P3\\ CANADA}

\email[Jean-Christophe Aval]{aval@math.u-bordeaux.fr}
\email[Nantel Bergeron]{bergeron@mathstat.yorku.ca}
\urladdr[Nantel Bergeron]{http://www.math.yorku.ca/bergeron}

\date{\today}
\thanks{N. Bergeron is supported in part by NSERC and PREA}


\begin{document}

\begin{abstract}
A lattice diagram is a finite list $L=\big((p_1,q_1),\ldots ,(p_n,q_n)\big)$ of
lattice cells. The corresponding lattice diagram determinant is
$\Delta_L(X;Y)=\det \|\, x_i^{p_j}y_i^{q_j}\, \|$. The
space $M_L$ is the space spanned by all partial derivatives of
$\Delta_L(X;Y)$. We describe here how a Schur function partial derivative
operator acts on lattice
diagrams with distinct cells in the positive quadrant.
 \end{abstract}

\maketitle

\section{Introduction}

The lattice cell in the $i+1^{st}$ row and $j+1^{st}$ column of the positive 
quadrant of the plane is denoted by $(i,j)$.
We order the set of all lattice cells using the following {\sl
lexicographic} order:
   \begin{equation}\label{lex}
    (p_1,q_1)<(p_2,q_2) \quad\iff\quad q_1<q_2\quad\hbox{or}
     \quad [q_1=q_2\hbox{ and }p_1<p_2].
   \end{equation}
For our purpose, a {\sl lattice diagram} is a finite list
$L=\big((p_1,q_1),\ldots
,(p_n,q_n)\big)$ of lattice cells such that $(p_1,q_1)\le (p_2,q_2)\le\ldots
\le(p_n,q_n)$. Following the definitions and conventions of
\cite{Berg et
al}, the coordinates $p_i$ and $q_i$ of a cell
$(p_i,q_i)$ indicate the row and
column  position,  respectively, of the cell.
For $\mu_1\geq \mu_2\geq \cdots \geq \mu_k>0$, we say that $\mu=(\mu_1,
\mu_2, \ldots,
\mu_k)$ is a {\sl partition} of $n$ if $n=\mu_1+\cdots +\mu_k$. We associate to
a partition $\mu$ the following lattice (Ferrers) diagram $\big((i,j)\,:\,0\le
i\le k-1,\,0\le j\le \mu_{i+1}-1\big)$, distinct cells ordered with
\ref{lex}, and we use the symbol
$\mu$ for both the partition and its associated Ferrers diagram.
For example, given the partition $(4,2,1)$, its Ferrers diagram is:
  $$\Young{\Carre{2,0}\cr
           \Carre{1,0}&\Carre{1,1}\cr
           \Carre{0,0}&\Carre{0,1}&\Carre{0,2}&\Carre{0,3}\cr
           }\quad.$$
This consists of the lattice cells
$\big((0,0),(1,0),(2,0),(0,1),(1,1),(0,2),(0,3)\big)$.

Given a lattice diagram $L=\big((p_1,q_1), (p_2,q_2),\ldots , (p_n,q_n)\big)$
we define the {\sl lattice diagram determinant}
  $$
  \Delta_L(X;Y)= \det \left\|
{{x_i^{p_j}y_i^{q_j}}\over{p_j!q_j!}}\right\|_{i,j=1}^n\,,
  $$
where $X=x_1,x_2,\ldots,x_n$ and $Y=y_1,y_2,\ldots,y_n$. This determinant
clearly vanishes if any
cell has multiplicity greater than one, and we set $\Delta_L(X;Y)=0$ if a
coordinate of any
cell is negavive. The determinant
$\Delta_L(X;Y)$ is  bihomogeneous of degree $|p|=p_1+\cdots +p_n$ in $X$
and  degree $|q|=q_1+\cdots +q_n$  in  $Y$. The factorials will
ensure that the lattice diagram determinants behave nicely under partial
derivatives.

For a polynomial $P(X;Y)$ we denote by $P(\partial X;\partial Y)$ the
differential
operator obtained from $P$,
substituting every variable $x_i$ by the operator
$\partial\,\,\,\,\over\partial x_i$ and every variable $y_j$ by the operator
$\partial\,\,\,\,\over\partial y_j$. A permutation $\sigma\in \S_n$ acts
diagonally on
a polynomial $P(X;Y)$ as follows:
$\sigma  P(X;Y)\,=\, P(x_{\sigma_1},x_{\sigma_2},\ldots
,x_{\sigma_n};y_{\sigma_1},y_{\sigma_2},\ldots ,y_{\sigma_n})$.
Under this action, $\Delta_L(X;Y)$ is clearly an alternant.

These lattice diagram determinants are crucial in the study of the so-called 
``$n!$ conjecture'' of A. Garsia and M. Haiman \cite{gh}, recently proven by 
M. Haiman \cite{geo}, and in generalizations of this question (see 
\cite{aval2,Berg et al} for example). To be more precise the key object is the
 vector space spanned by all partial derivatives of a given lattice diagram 
determinant $\Delta_L$, which we denote by 
$${\bf M}_L={\mathcal L}_\partial[\Delta_L].$$ 

Very useful in the comprehension of the structure of the ${\bf M}_L$ spaces 
are the ``shift operators''. These operators are special symmetric derivative 
operators, whose action on the lattice diagram determinants could be easily 
described in terms of movements of cells.

Another interest related to the shift operators is the hope to obtain a 
description of the vanishing ideal of ${\bf M}_L$, which is defined as:
$${\mathcal I}_L=\{f\in\Q[X;Y] ;\ f(\partial X;\partial Y)\Delta_L(X;Y)=0\}.$$
The structure of ${\bf M}_L$ and of ${\mathcal I}_L$ are closely related and 
the shift operators are crucial tools to study ${\mathcal I}_L$ (see 
\cite{a3,aval2,a7} for some applications).

Let us recall results of
\cite{aval2} that describe the effects of power sums, elementary and 
homogeneous symmetric
differential operators on
lattice diagram determinants. 

For the sake of simplicity, we limit our descriptions to $X$-operators; the
$Y$-operators are
similar. Recall that
  $$\displaystyle P_k(X)=\sum_{i=1}^k x_{i}^k$$
$$\displaystyle e_k(X)=\sum_{1\le i_1 < i_2 <\cdots<i_k\le n} x_{i_1}
    x_{i_2}\cdots x_{i_k}$$
$$\displaystyle h_k(X)=\sum_{1\le i_1 \le i_2 \le\cdots\le i_k\le n} x_{i_1}
    x_{i_2}\cdots x_{i_k}$$
are respectively the $k$-th power sum, elementary and homogeneous symmetric 
polynomial. 

Now, to state the next proposition, we need to introduce some notation. For a 
lattice diagram $L$, we denote by $\overline L$ its complement in the
positive quadrant
(it is an infinite subset). Again we order $\overline L=\{(\overline
p_1,\overline q_1), (\overline
p_2,\overline q_2),\dots\,\}$ using the lexicographic order \ref{lex}.
Let
$L$ be a lattice diagram with $n$ distinct cells in the positive quadrant. For
any integer $k\ge 1$ we have:

\begin{prop}[Propostion I.1 \cite{Berg et al}, Propositions 2.4, 2.6 
\cite{aval2}]\label{elem}

\begin{equation}\label{pkop}
P_k(\partial X)\Delta_L(X,Y)=\sum_{i=1}^n\pm\Delta_{P_k(i;L)}(X,Y),
\end{equation}
where ${P_k(i;L)}$ is the diagram obtained by replacing the $i$-th biexponent 
$(p_i,q_i)$ by $(p_i-k,q_i)$ and the coefficient $\epsilon(L,P_k(i;L))$ is a 
positive integer. The sign in \ref{pkop} is the sign of the permutation that 
reorders the obtained biexponents with respect to the lexicographic order 
\ref{lex}.

\begin{equation}
e_k(\partial X)\Delta_L(X;Y)=\sum_{1\le i_1<i_2<\cdots<i_k\le
n}\Delta_{e_k(i_1,\ldots,i_k;L)}(X;Y)
\end{equation}
 where
$e_k(i_1,\ldots,i_k;L)$ is the lattice diagram obtained from $L$ by
replacing the biexponents
$(p_{i_1},q_{i_1}),\ldots,(p_{i_k},q_{i_k})$ with
$(p_{i_1}-1,q_{i_1}),\ldots,(p_{i_k}-1,q_{i_k})$.

\begin{equation}\label{eqh}
h_k(\partial X)\Delta_L(X,Y)=\sum_{1\le i_1<i_2<\cdots<i_k}
  \Delta_{h_k(i_1,\dots,i_k;L)}(X,Y)
\end{equation}
where
$h_k(i_1,\dots,i_k;L)$ is the lattice diagram with the following
complement diagram. Replace the
biexponents
$(\overline p_{i_1},\overline q_{i_1}),\dots,$ $(\overline
p_{i_k},\overline q_{i_k})$ of the complement
$\overline L$ with
$(\overline p_{i_1}+1,\overline q_{i_1}),\\ \dots,(\overline
p_{i_k}+1,\overline q_{i_k})$ and keep the other unchanged. 
\end{prop}

The aim of this work is to obtain a description similar to the previous 
proposition of the effect of a partial Schur differential symmetric operator 
on a lattice diagram determinant. We obtain such a result in the next section 
and prove it.

\section{Schur operators}

 Following \cite{Mac}, recall that for a partition
$\lambda=(\lambda_1,\lambda_2,\ldots,\lambda_k)$ the conjugate (transpose)
partition is denoted by
$\lambda'=(\lambda'_1,\lambda'_2,\ldots,\lambda'_\ell)$. With this in mind,
the Schur polynomial
indexed by
$\lambda$ is

$$S_\lambda(X)=\det{\|\,e_{\lambda'_j+i-j}(X)
\,\|}
$$
with the understanding that $e_0(X)=1$ and $e_k(X)=0$ if
$k<0$. The Schur
polynomials also have a description in terms of column-strict Young tableaux.
Given $\lambda$ a
partition of $n$, a tableau of shape $\lambda$ is a map
$T\colon\lambda\to\{1,2,\ldots,n\}$. We say that $T$ is a column-strict
Young tableau if it
is weakly increasing
along the rows and strictly increasing along the columns of $\lambda$. That is
$T(i,j)\le T(i,j+1)$ and
$T(i,j)<T(i+1,j)$ wherever it is defined. We denote by ${\mathcal
T}_\lambda$ the set of all
column-strict Young tableaux of shape $\lambda$. For any tableau $T$, we define
$X^T=\prod_{i=1}^n
x_i^{|T^{-1}(i)|}$. As seen in~\cite{Mac}, we have
  $$S_\lambda(X)=\sum_{T\in{\mathcal T}_\lambda} X^T.
$$

It is convenient to define the following function on lattice diagrams.
\begin{equation}\label{coef}\epsilon(L) = \left\{
     \begin{array}{rll}
        1&& \hbox{if $L$ has $n$ distinct cells in the positive quadrant,}\\
      &&\\
        0&& \hbox{otherwise.}
      \end{array}
     \right.
\end{equation}

Let  $L$ be a lattice diagram with $n$ distinct cells in the positive
quadrant. For any partition
$\lambda$ of an integer $k\ge 1$ we have

\begin{thm}\label{schur}
$$S_\lambda(\partial X)\Delta_L(X;Y)=\sum_{T\in{\mathcal T}_\lambda}
  \epsilon'(T,L)\Delta_{\partial T(L)}(X;Y)$$
 where
$\partial T(L)$ is the lattice diagram obtained from $L$ by
replacing the biexponents
$(p_i,q_i)$ with
$(p_i-|T^{-1}(i)|,q_i)$ for $1\le i\le n$. The coefficient $\epsilon'(T,L)$
is described 
as follows. Let $T_1,T_2,\ldots,T_\ell$ be the $\ell$ columns of $T$ then
$\partial T(L)=\partial
T_1\partial T_2\cdots\partial T_\ell(L)$ and
 \begin{equation}\label{scoef}\epsilon'(T,L)=
    \epsilon\big(\partial T(L)\big)\ \cdots\
    \epsilon\big(\partial T_{\ell-1}\partial T_\ell(L)\big)\
\epsilon\big(\partial T_\ell(L)\big)
 \end{equation}
where $\epsilon$ is defined in \ref{coef}. Hence $\epsilon'(T,L)$ is 0 or 1.
\end{thm}

We shall prove this result using the Proposition \ref{elem} and an
adaptation of the involution
defined in \cite{RS}. We will see in the proof at the end of this section
that the succesive order
in wich we apply the operators
$\partial T_j$ to the lattice diagram $L$ in the equation \ref{scoef} is
not arbitrary. The
result and the proof depend on that precise order and this does not appear
in any previous work.

To start, we remark that the Theorem
\ref{schur} and  Proposition \ref{elem} agree on their domain of
definition. This is because
$e_k=S_{1^k}$ and the tableau of shape $1^k$ corresponds to sequences $1\le
i_1<i_2<\cdots<i_k\le
n$.  Now let
$\ell$ be the number of components of
$\lambda'$ and expand the determinant
 \begin{equation}\label{sexp}
   S_\lambda(X)=\det{\|\,e_{\lambda'_j+i-j}(X)\,\|} = \sum_{\sigma\in
\S_\ell } sgn(\sigma)
   e_{\sigma(\lambda'+\delta_\ell)-\delta_\ell}.
 \end{equation}
Here $\delta_\ell=(\ell-1,\ell-2,...,1,0)$ and $e_{\alpha_j}=0$ if
$\alpha_j<0$. If we have
$\alpha=\alpha_1,\alpha_2,\ldots,\alpha_\ell$ a sequence of integer we let
$e_\alpha=e_{\alpha_1}e_{\alpha_{2}}\cdots e_{\alpha_\ell}$. Here the order
in which we write
this product matters. For $\ell=1$, as noted before,  Proposition
\ref{elem} can be rewritten as
\begin{equation}\label{eq:1} e_{\alpha_1}(\partial
X)\Delta_L(X;Y)=\sum_{T_1\in\,{\mathcal
    T}_{1^{^{\alpha_1}}}} \epsilon(\partial T_1(L))\Delta_{\partial
T_1(L)}(X;Y).
\end{equation}
Where ${\mathcal T}_{1^{^{\alpha_1}}}$ is the set of $\alpha_1$-column
tableaux   with content
in $\{1,2,\ldots,n\}$, strictly increasing in the column. Here $
\epsilon'(T_1,L)=
\epsilon(\partial T_1(L))$. Suppose now that $\ell=2$. We use \ref{eq:1} with
$e_{\alpha_2}(\partial X)$ and apply $e_{\alpha_1}(\partial X)$ on both
side. That gives
 \begin{eqnarray*}
  e_{\alpha}(\partial X)\Delta_L(X;Y)&=&e_{\alpha_1}(\partial X)
e_{\alpha_2}(\partial
X)\Delta_L(X;Y) \\
 &=& \sum_{T_2\in\,{\mathcal T}_{1^{^{\alpha_2}}}} \epsilon(\partial
T_2(L)) e_{\alpha_1}(\partial
X)
\Delta_{\partial T_2(L)}(X;Y)\\
 &=& \sum_{T_1\in\,{\mathcal T}_{1^{^{\alpha_1}}}} \sum_{T_2\in\,{\mathcal
T}_{1^{^{\alpha_2}}}}
\epsilon(\partial T_2(L))\epsilon(\partial T_1\partial T_2(L))
\Delta_{\partial T_1\partial T_2(L)}(X;Y)\\
 \end{eqnarray*}
Now let
$\mathcal{CT}_{\alpha}=\mathcal{CT}_{\alpha_1,\alpha_2,\ldots,\alpha_\ell}$
be the set of $\ell$
columns ${\bf T}=(T_1,T_2,\ldots,T_\ell)$ where $T_j\in\,{\mathcal
T}_{1^{^{\alpha_j}}}$. We can
represent $\bf T$ as a tableau $\alpha\to\{1,2,\ldots,n\}$ where as before
we identify the
composition
$\alpha$ with the lattice diagram
$\big((i,j)\,|\,0\le i\le\alpha_{j+1}-1,\ 0\le j\le\ell-1\big)$, with
distinct cells ordered by
\ref{lex}. The tableau
$\bf T$ is strictly increasing along every columns and has no restriction
along rows.
Note that the shape
$\alpha$ is not necessarily a partition. We can now simplify our
computation above  and write for
$\ell=2$:
 \begin{equation}\label{ncelem}
   e_{\alpha}(\partial X)\Delta_L(X;Y) = \sum_{{\bf T}\in
\mathcal{CT}_{\alpha}} \epsilon'({\bf
T},L) \Delta_{\partial {\bf T}(L)}(X;Y),
 \end{equation}
where $\partial {\bf T}(L)=\partial T_1\partial T_2\cdots\partial
T_\ell(L)$ is the lattice diagram
obtained from $L$ by replacing the biexponents $(p_i,q_i)$ with
$(p_i-|{\bf T}^{-1}(i)|,q_i)$ for $1\le i\le n$ and
 \begin{equation}\label{tcoef}
  \epsilon'({\bf T},L)=    \epsilon\big(\partial {\bf T}(L)\big)\ \cdots\
    \epsilon\big(\partial T_{\ell-1}\partial T_\ell(L)\big)\
\epsilon\big(\partial
    T_\ell(L)\big).
 \end{equation}
It is clear, by induction, that this is true for all $\ell\ge 2$
as well. We must also remark here that if one of the $\alpha_j<0$ the sum
\ref{ncelem}
must be set to zero.

We can now start the computation of the operator \ref{sexp} using \ref{ncelem}:
 \begin{eqnarray}
  S_\lambda(\partial X)\Delta_L(X;Y) &=& \sum_{\sigma\in\S_\ell}
sgn(\sigma) e_{\sigma(\lambda'+\delta_\ell)-\delta_\ell}(\partial X)
\Delta_{L}(X;Y)\nonumber\\
 &=& \sum_{\sigma\in\S_\ell}\quad \sum_{{\bf T}\in
\mathcal{CT}_{\sigma(\lambda'+\delta_\ell)-\delta_\ell}} sgn(\sigma)
 \epsilon'({\bf T},L) \Delta_{\partial {\bf T}(L)}(X;Y).\label{cont}
 \end{eqnarray}

Now we need to construct an involution on the set indexing the double sum
such that all term cancels,
unless $\sigma$ is the identity and $\bf T\in\mathcal{T}_\lambda$. Here is
an example of a ${\bf
T}\in{\mathcal{CT}}_{(1,0,3,2,4,1)}$
  $${\bf T}\quad=\quad\Young{
           &&&&\Carre{8}\cr
           &&\Carre{10}&&\Carre{6}\cr
           &&\Carre{8}&\Carre{9}&\Carre{5}\cr
           \Carre{2}&\Vide{}&\Carre{7}&\Carre{3}&\Carre{4}&\Carre{4}\cr}
  $$
The only requirement is that $\bf T$ is strictly increasing in columns.

Let us first concentrate on $\ell=2$ and let
$\lambda'=(\lambda'_1,\lambda'_2)$. We have two
possible shapes $\alpha=(\alpha_1,\alpha_2)$, either $\lambda'=
Id(\lambda'+\delta_2)-\delta_2$ or
$(\lambda'_2-1,\lambda'_1+1)=(1,2)(\lambda'+\delta_2)-\delta_2$, where $(i,j)$ 
is the usual notation for transpositions. These two cases are
completely characterized
by $\alpha_1<\alpha_2$ or $\alpha_1\ge\alpha_2$. We now define an
involution similar to \cite{RS}.

Given $\bf{T} \in \mathcal{CT}_\alpha$ we associate two words $w_{_{\bf
T}}$ and $\widehat{w}_{_{\bf
T}}$. This method is originally due to A. Lascoux and M.-P. Sch\"utzenberger 
(cf. \cite{las}). The first $w_{_{\bf T}}$ is all the entries ${\bf T}(i,j)$ 
of ${\bf T}$ sorted in increasing
order. For example if
  $${\bf T}\quad=\quad\Young{
           &\Carre{9}\cr
           &\Carre{6}\cr
           \Carre{9}&\Carre{5}\cr
           \Carre{3}&\Carre{4}\cr}
  $$
then $w_{_{\bf T}}= 3\,4\,5\,6\,9\,9$. Now we associate to $w_{_{\bf T}}$
its parentheses
structure $\widehat{w}_{_{\bf T}}$. For this, we list the entries in
$w_{_{\bf T}}$ and associate to
an entry from the first column of $\bf T$ a left parenthesis, and to an
entry of the second column a
right parenthesis. For two columns, the same entry appears at most twice,
in which case the first
one we read in
$w_{_{\bf T}}$ is assumed to be from the first column of $\bf T$. In the
example above
$w_{_{\bf T}}= 3\,4\,5\,6\,9\,9$ and $\widehat{w}_{_{\bf
T}}=(\,)\,)\,)\,(\,)$.

There is a natural way to pair parentheses under the usual rule of
parenthesization. In any word
$\widehat{w}_{_{\bf T}}$ some parentheses will be paired and other will be
unpaired. In our example,
 $\widehat{w}_{_{\bf T}}=(\,)\, \big)\,\big)\,(\,)$, the first two
parentheses and last two are
paired and the two parentheses in the midle are unpaired. The subword of
any $\widehat{w}_{_{\bf
T}}$ consisting of unpair parentheses must be of the form
$)\,)\cdots)\,(\cdots (\,($

We the have the following useful result.
\begin{prop}\cite[Proposition 5]{RS}
\label{shimo}
 A tableau ${\bf T}=(T_1,T_2,\ldots,T_\ell)\in\mathcal{CT}_\alpha$ is a
column strict Young tableau
${\bf T}\in\mathcal{T}_\lambda$ if and only if there are no unpaired right
parentheses in
$\widehat{w}_{_{T_j,T_{j+1}}}$ for all $1\le j\le \ell-1$ and two columns
$T_j,T_{j+1}$ of $\bf T$.
\end{prop}

Remark here that if $\alpha=(\alpha_1,\alpha_2,\ldots,\alpha_\ell)$ is not
a partition, that is
$\alpha_j<\alpha_{j+1}$ for some $1\le j\le\ell-1$, then nescessarily
$\widehat{w}_{_{T_j,T_{j+1}}}$
will contain more right parentheses than left parentheses and some will be
left unpaired and no
$\bf T\in\mathcal{CT}_\alpha$ could be a column strict Young tableau.

We return to the construction of the involution from to \cite{RS} for
$\lambda'=(\lambda'_1,\lambda'_2)$. Let

$$A=\mathcal{CT}_{(\lambda'_1,\lambda'_2)}\cup\mathcal{CT}_{(\lambda'_2-1,
\lambda'_1+1)}.$$
The involution is a map $\Psi\colon A\to A$ defined as follows. Let ${\bf
T}\in\mathcal{CT}_{(\alpha_1,\alpha_2)}\subset A$ and consider
$\widehat{w}_{_{\bf T}}$. The subword of unpaired parentheses contains
$r\ge 0$ unpaired right
parentheses followed by $l\ge 0$ unpaired left parentheses. We have that
$l-r=\alpha_1-\alpha_2$.

\noindent $\bullet$ If $r=0$, then ${\bf
T}\in\mathcal{T}_\lambda\subset \mathcal{CT}_{\lambda'}$ and we define
$\Psi({\bf T})={\bf T}$.

\noindent $\bullet$  If $l\ge r>0$, then ${\bf
T}\in\mathcal{CT}_{\lambda'}\setminus
\mathcal{T}_\lambda$ and we define $\Psi({\bf T})={\bf
T'}\in\mathcal{CT}_{(\lambda'_2-1,\lambda'_1+1)}$, the unique tableau such that
${w}_{_{\bf T'}}={w}_{_{\bf T}}$ and $\widehat{w}_{_{\bf T'}}$ is obtained
from $\widehat{w}_{_{\bf
T}}$ replacing the $l-r+1$ leftmost unpaired left parentheses by right
parentheses.

\noindent $\bullet$  If $r>l$, then ${\bf
T}\in\mathcal{CT}_{(\lambda'_2-1,\lambda'_1+1)}$ and we
define $\Psi({\bf T})={\bf T'}\in\mathcal{CT}_{\lambda'}\setminus
\mathcal{T}_\lambda$, the unique tableau such that
${w}_{_{\bf T'}}={w}_{_{\bf T}}$ and $\widehat{w}_{_{\bf T'}}$ is obtained
from $\widehat{w}_{_{\bf
T}}$ replacing the $r-l-1$ rightmost unpaired right parentheses by left
parentheses.

Now in the general case, that is if $\ell\ge 2$, let
  $$A=\bigcup_{\sigma\in\S_\ell}
\mathcal{CT}_{\sigma(\lambda'+\delta_\ell)-\delta_\ell}.$$
For ${\bf T}\in\mathcal{CT}_\alpha\subset A$, the composition $\alpha$
completely characterizes the
permutation $\sigma\in\S_\ell$ such that
$\alpha=\sigma(\lambda'+\delta_\ell)-\delta_\ell$. In
particular $\alpha$ is a partition if and only if $\sigma=Id$. We read the
rows of $\bf T$ from
right to left, bottom to top. We find this way the first  pair $(i,j)$ and
$(i,j+1)$ such that
$$\hbox{$T(i,j)>T(i,j+1)$\quad or\qquad  $(i,j)\not\in \alpha$ and
$(i,j+1)\in
\alpha$.}$$

\noindent $\bullet$ If there is no such pair, then we have ${\bf
T}\in\mathcal{T}_\lambda\subset
\mathcal{CT}_{\lambda'}$ and we define $\Psi({\bf T})={\bf T}$.

\noindent $\bullet$ If we find such a pair, then we have ${\bf
T}\in\mathcal{CT}_\alpha\subset A\setminus
\mathcal{T}_\lambda$. We
define $\Psi({\bf T})={\bf T'}\in\mathcal{CT}_{\beta}\subset A\setminus
\mathcal{T}_\lambda$ where $\bf T'$ is obtained from $\bf T$ using the
procedure above
to the two columns
$T_{j+1},T_{j+2}$. By construction if
$\alpha=\sigma(\lambda'+\delta_\ell)-\delta_\ell$, then
$\beta=\sigma(j,j+1)(\lambda'+\delta_\ell)-\delta_\ell$.

The fact that $\Psi$ is a well defined involution is done in several papers,
for example, in \cite{RS}, section 3. Let us give one example. For
  $${\bf T}\quad=\quad\Young{
           &&\Carre{8}\cr
           \Carre{10}&&\Carre{6}\cr
           \Carre{8}&\Carre{9}&\Carre{5}\cr
           \Carre{7}&\Carre{3}&\Carre{4}\cr}
   \hbox{ \quad we have \quad}
     \Psi({\bf T})\quad=\quad\Young{
           \Carre{10}&\Carre{9}&\Carre{6}\cr
           \Carre{8}&\Carre{8}&\Carre{5}\cr
           \Carre{7}&\Carre{3}&\Carre{4}\cr
           }.
  $$
The pair $(1,2)$ and $(1,3)$ is the first one where $T(1,2)>T(1,3)$. We
thus apply the involution
on the second and third column. We have here
${w}_{_{T_2,T_3}}=3\,4\,5\,6\,8\,9$ and
$\widehat{w}_{_{T_2,T_3}}=(\,)\,)\,)\,)\,($. There are $r=3$ unpaired right
parentheses followed by
$l=1$ unpaired left parenthesis. We must change $r-l-1=1$ unpaired left
parenthesis for a right one.
That is $\widehat{w}_{_{T'_2,T'_3}}=(\,)\,)\,)\,(\,($. That moves the entry 
$8$ from the third column to the second column.

\proof[Proof of Theorem \ref{schur}:]
We return to the computation \ref{cont} using the notation we have developped:
 \begin{eqnarray*}
  S_\lambda(\partial X)\Delta_L(X;Y)
 &=& \sum_{{\bf T}\in\,
\mathcal{CT}_{\sigma(\lambda'+\delta_\ell)-\delta_\ell}\subset A} sgn(\sigma)
 \epsilon'({\bf T},L) \Delta_{\partial {\bf T}(L)}(X;Y).
 \end{eqnarray*}
The involution constructed above matches the term in the sum corresponding
to  ${\bf
T}\in\mathcal{CT}_{\sigma(\lambda'+\delta_\ell)-\delta_\ell}\subset A\setminus
\mathcal{T}_\lambda$ with
${\bf
T'}\in\mathcal{CT}_{\sigma(j,j+1)(\lambda'+\delta_\ell)-\delta_\ell}\subset
A\setminus
\mathcal{T}_\lambda$.
Clearly, we have that $sgn(\sigma)=-sgn(\sigma(j,j+1))$ and $\partial {\bf
T}(L)=\partial {\bf
T'}(L)$. Once we show that
   \begin{equation}\label{invc}\epsilon'({\bf T},L)=\epsilon'({\bf
T'},L)\end{equation}
the Theorem \ref{schur} will follow from the fact that all the terms in $A
\setminus
\mathcal{T}_\lambda$ will cancel out and the remaining terms are in
$\mathcal{T}_\lambda$ with the
desired coefficient.

To establish \ref{invc} we need to show that if $\epsilon'({\bf T},L)\ne 0$
then $\epsilon'({\bf
T'},L)\ne 0$, for they will then both be equal to 1. From \ref{tcoef}
 $$  \epsilon'({\bf T},L)=   \epsilon'((T_1,T_2,\ldots,T_\ell),L)=
\epsilon(\partial{\bf T}(L))\
     \cdots\
    \epsilon(\partial T_{\ell -1}\partial T_\ell(L))\ \epsilon(\partial
T_\ell(L)).$$
Similarly $\epsilon'({\bf T'},L)=
\epsilon'((T_1,T_2,\ldots,T'_{j+1},T'_{j+2},\ldots,T_\ell),L)$
for some $0\le i\le\ell-1$. If  $\epsilon'({\bf T},L)\ne 0$, then
$\epsilon(\partial
T_k\cdots\partial T_\ell(L))=1$ for
$1\le k\le\ell$. For $1\le k\le j+1$ we clearly have:
  $$\epsilon(\partial T_k\cdots\partial T_{j+1}\partial
T_{j+2}\cdots\partial T_\ell(L))
  =\epsilon(\partial T_k\cdots\partial T'_{j+1}\partial
T'_{j+2}\cdots\partial T_\ell(L)).
 $$
For $j+3\le k\le \ell$, the corresponding terms of $\epsilon'({\bf T},L)$
and $\epsilon'({\bf
T'},L)$ are the same. Let $\tilde L=\partial T_{j+3}\cdots\partial
T_{\ell-1}\partial T_\ell L$, the
equality
\ref{invc} will follow as soon as we show that
   \begin{equation}\label{invc2}\epsilon(\partial T_{j+2}(\tilde
L))=1\quad\hbox{and}\quad
\epsilon(\partial T_{j+1}\partial T_{j+2}(\tilde L))=1\qquad\implies\qquad
\epsilon(\partial
T'_{j+2}(\tilde L))=1\end{equation}
for all $\tilde L$ such that $\epsilon(\tilde L)=1$.

Let $\beta=\sigma(j,j+1)(\lambda'+\delta_\ell)-\delta_\ell$, the shape of
$\bf T$. Suppose that
$\epsilon(\partial T'_{j+2}(\tilde L))=0$. This implies that there is an entry
$1\le k={\bf T'}(i,j+2)\le n$ such that the cells $(p_k,q_k)\in\tilde L$ and
$(p_{k-1},q_{k-1})=(p_k-1,q_k)\in \tilde L$, and $k-1\ne{\bf T'}(i-1,j+2)$
is not an entry of
$T'_{j+2}$. Now since $\epsilon(\partial T_{j+1}\partial T_{j+2}(\tilde
L))=1$ we must have that
both $k$ and $k-1$ are entries of $T_{j+1},T_{j+2}$. This implies that
$k-1$ is an entry of
$T'_{j+1}$ and $k$ is not. This analysis shows that $k-1$ and $k$ are
entries of
$w_{_{T'_{j+1}T'_{j+2}}}$ with multiplicity one,  $k-1$ is in the column
$T'_{j+1}$ and $k$ is
in the column $T'_{j+2}$. They will be consecutive entries in
$w_{_{T'_{j+1}T'_{j+2}}}$ and will be
paired in $\widehat{w}_{_{T'_{j+1}T'_{j+2}}}$. This would imply that
$T_{j+2}$ in $\Psi({\bf
T'})={\bf T}$ contains the entry $k$ but not $k-1$ and $\epsilon(\partial
T_{j+2}(\tilde L))=0$,
contrary to our hypothesis. This completes our proof.\endproof

\begin{rem}\rm Given a lattice diagram $L$ and a column strict tableau
$T\in\mathcal{T}_\lambda$, we have that $\epsilon'(T,L)=1$ exactly when we
can {\sl move} the cells
of $L$ by one, reading $T$ column by column, from right to left, without
having any cells
colliding.
\end{rem}

\begin{cor} For $h_k(X)=s_{(k)}(X)$ we have
 $$h_k(\partial X)\Delta_L(X;Y)=\sum_{1\le j_1\le j_2\le\cdots\le j_k\le
n}\epsilon'((j_1,\ldots,j_k),L)\Delta_{\partial_{j_1}\cdots\partial_{j_k}(L)}
(X;Y)$$
\end{cor}

This is equivalent to the description in \cite{aval2}. The only way that
$\epsilon'((j_1,\ldots,j_k),L)\ne 0$ is if the cells $j_1,\ldots,j_k$ that
moves down are moved
into holes. This can be described as holes moving up.

\end{document}